\documentclass[reqno,a4paper]{amsart}
\usepackage{amssymb,ifpdf}
\ifpdf
 \usepackage[hyperindex,pagebackref]{hyperref}%
\else
 \expandafter\ifx\csname dvipdfm\endcsname\relax
 \usepackage[hypertex,hyperindex,pagebackref]{hyperref}
 \else
 \usepackage[dvipdfm,hyperindex,pagebackref]{hyperref}
 \fi
\fi
\newtheorem{thm}{Theorem}[section]

\theoremstyle{remark}
\newtheorem{rem}{Remark}[section]
\allowdisplaybreaks[4]
\DeclareMathOperator{\td}{d}
\numberwithin{equation}{section}

\begin{document}

\title[A formula for Bernoulli numbers of the second kind]
{An explicit formula for computing Bernoulli numbers of the second kind in terms of Stirling numbers of the first kind}

\author[F. Qi]{Feng Qi}
\address[F. Qi]{Department of Mathematics, College of Science, Tianjin Polytechnic University, Tianjin City, 300387, China; Institute of Mathematics, Henan Polytechnic University, Jiaozuo City, Henan Province, 454010, China}
\email{\href{mailto: F. Qi <qifeng618@gmail.com>}{qifeng618@gmail.com}, \href{mailto: F. Qi <qifeng618@hotmail.com>}{qifeng618@hotmail.com}, \href{mailto: F. Qi <qifeng618@qq.com>}{qifeng618@qq.com}}
\urladdr{\url{http://qifeng618.wordpress.com}}

\begin{abstract}
In the paper, the author finds an explicit formula for computing Bernoulli numbers of the second kind in terms of Stirling numbers of the first kind.
\end{abstract}



\thanks{This paper was typeset using \AmS-\LaTeX}

\maketitle

\section{Introduction}

It is well known that Stirling numbers of the first kind $s(n,k)$ for $n\ge k\ge1$ may be generated by
\begin{equation}\label{s(n-m)-1dfn}
\frac{[\ln(1+x)]^k}{k!}=\sum_{n=k}^\infty s(n,k)\frac{x^n}{n!},\quad |x|<1
\end{equation}
and that Bernoulli numbers of the second kind $b_n$ for $n\ge0$ may be generated by
\begin{equation}\label{bernoulli-second-dfn}
\frac{x}{\ln(1+x)}=\sum_{n=0}^\infty b_nx^n.
\end{equation}
\par
In~\cite{Nemes-JIS-2011}, the following formula for computing Bernoulli numbers of the second kind in terms of Stirling numbers of the first kind was derived:
\begin{equation}
b_n=\frac1{n!}\sum_{k=0}^n\frac{s(n,k)}{k+1}.
\end{equation}
\par
The main aim of this paper is to find a new and explicit formula for computing Bernoulli numbers of the second kind in terms of Stirling numbers of the first kind. The main result of this paper may be stated as the following theorem.

\begin{thm}\label{Bern-2rd1stStirl=thm}
For $n\ge2$, Bernoulli numbers of the second kind $b_n$ may be computed in terms of Stirling numbers of the first kind $s(n,k)$ by
\begin{equation}\label{Bern-2rd=formulaStirl}
b_n=\frac1{n!}\sum_{k=1}^{n-1}(-1)^{k} \frac{s(n-1,k)}{(k+1)(k+2)}.
\end{equation}
\end{thm}

As a remark, a relation between the harmonic numbers and Stirling numbers of the first kind $s(n,2)$ is also derived.

\section{Proof of Theorem~\ref{Bern-2rd1stStirl=thm}}

The proof of Theorem~\ref{Bern-2rd1stStirl=thm} is based on some results elementarily and inductively obtained in~\cite{Filomat-36-73-1.tex} and its preprint~\cite{2rd-Bernoulli-2012.tex}. These results can be recited as follows.
\begin{enumerate}
\item
Corollary~2.3 in~\cite{Filomat-36-73-1.tex, 2rd-Bernoulli-2012.tex} states that Stirling numbers of the first kind $s(n,k)$ for $1\le k\le n$ may be computed by
\begin{equation}\label{s(n,k)-sum}
s(n,k)=(-1)^{n+k}(n-1)!\sum_{\ell_1=1}^{n-1} \frac1{\ell_1}\sum_{\ell_2=1}^{\ell_1-1}\frac1{\ell_2}\dotsm \sum_{\ell_{k-2}=1}^{\ell_{k-3}-1}\frac1{\ell_{k-2}} \sum_{\ell_{k-1}=1}^{\ell_{k-2}-1}\frac1{\ell_{k-1}}.
\end{equation}
This formula may be reformulated as
\begin{equation}
(-1)^{n-k}\frac{s(n,k)}{(n-1)!}= \sum _{m=k-1}^{n-1}\frac1m\biggl[(-1)^{m-(k-1)}\frac{s(m,k-1)}{(m-1)!}\biggr].
\end{equation}
\item
Corollary~2.4 in~\cite{Filomat-36-73-1.tex, 2rd-Bernoulli-2012.tex} reads that for $1\le k\le n$ Stirling numbers of the first kind $s(n,k)$ satisfies the recursion
\begin{equation}\label{1st-Stirling-recurs}
s(n+1,k)=s(n,k-1) -ns(n,k).
\end{equation}
This is a recovery of the triangular relation for $s(n,k)$.
\item
Theorem~3.1 in~\cite{Filomat-36-73-1.tex, 2rd-Bernoulli-2012.tex} tells that Bernoulli numbers of the second kind $b_n$ for $n\ge2$ may be computed by
\begin{equation}\label{Bernulli-2rd=formula}
b_n=(-1)^n\frac1{n!}\Biggl(\frac1{n+1}+\sum_{k=2}^{n} \frac{a_{n,k}-na_{n-1,k}}{k!}\Biggr),
\end{equation}
where
\begin{equation}\label{a(n,2)-eq}
a_{n,2}=(n-1)!
\end{equation}
and, for $n+1\ge k\ge3$,
\begin{equation}\label{a=n=i=eq}
a_{n,k}=(k-1)!(n-1)!\sum_{\ell_1=1}^{n-1} \frac1{\ell_1}\sum_{\ell_2=1}^{\ell_1-1}\frac1{\ell_2}\dotsm \sum_{\ell_{k-3}=1}^{\ell_{k-4}-1}\frac1{\ell_{k-3}} \sum_{\ell_{k-2}=1}^{\ell_{k-3}-1}\frac1{\ell_{k-2}}.
\end{equation}
\end{enumerate}
\par
Observing the expressions~\eqref{s(n,k)-sum} and~\eqref{a=n=i=eq}, we obtain
\begin{equation}\label{a-stirl-relation}
a_{n,k}=(-1)^{n+k-1}(k-1)!s(n,k-1),\quad n+1\ge k\ge2.
\end{equation}
See~\cite[(2.18)]{Filomat-36-73-1.tex} and~\cite[(6.7)]{1st-Sirling-Number-2012.tex}. By this and the recursion~\eqref{1st-Stirling-recurs}, it follows that
\begin{align*}
a_{n,k}-na_{n-1,k}&=(-1)^{n+k-1}(k-1)![s(n,k-1)+ns(n-1,k-1)]\\
&=(-1)^{n+k-1}(k-1)![s(n-1,k-1)+s(n-1,k-2)].
\end{align*}
Substituting this into~\eqref{Bernulli-2rd=formula} reveals that
\begin{align*}
b_n&=\frac{(-1)^n}{n!}\Biggl(\frac1{n+1}+\sum_{k=2}^{n} \frac{(-1)^{n+k-1}[s(n-1,k-1)+s(n-1,k-2)]}{k}\Biggr)\\
&=\frac{(-1)^n}{(n+1)!}+\frac1{n!}\sum_{k=2}^{n} \frac{(-1)^{k-1}[s(n-1,k-1)+s(n-1,k-2)]}{k}\\
&=\frac{(-1)^n}{(n+1)!}+\frac1{n!} \Biggl[\sum_{k=2}^{n} \frac{(-1)^{k-1}s(n-1,k-1)}{k} +\sum_{k=2}^{n} \frac{(-1)^{k-1}s(n-1,k-2)}{k}\Biggr]\\
&=\frac{(-1)^n}{(n+1)!}+\frac1{n!} \Biggl[\sum_{k=2}^{n} \frac{(-1)^{k-1}s(n-1,k-1)}{k} +\sum_{k=1}^{n-1} \frac{(-1)^{k}s(n-1,k-1)}{k+1}\Biggr]\\
&=\frac{(-1)^n}{(n+1)!}+\frac1{n!}\frac{(-1)^{n-1}}n+\frac1{n!} \sum_{k=2}^{n-1} (-1)^{k-1}s(n-1,k-1)\biggl(\frac1k- \frac1{k+1}\biggr)\\
&=\frac1{n!}(-1)^{n-1}\biggl(\frac1n-\frac1{n+1}\biggr)+\frac1{n!} \sum_{k=2}^{n-1} (-1)^{k-1}s(n-1,k-1)\biggl(\frac1k- \frac1{k+1}\biggr)\\
&=\frac1{n!} \sum_{k=2}^{n} (-1)^{k-1}s(n-1,k-1)\biggl(\frac1k- \frac1{k+1}\biggr)\\
&=\frac1{n!}\sum_{k=2}^{n}(-1)^{k-1} \frac{s(n-1,k-1)}{k(k+1)}.
\end{align*}
Notice that in the above argument, we use the convention $s(n,0)=0$ for $n\in\mathbb{N}$ and the fact $s(n,n)=1$ for $n\ge0$. The proof of Theorem~\ref{Bern-2rd1stStirl=thm} is complete.

\section{Remarks}

In this section, we show some new findings by several remarks.

\begin{rem}
The idea in Theorem~\ref{Bern-2rd1stStirl=thm} and its proof ever implicitly thrilled through in~\cite[Remark~6.7]{1st-Sirling-Number-2012.tex}.
\end{rem}

\begin{rem}
Making use of the relation~\eqref{a-stirl-relation} in~\cite[Theorem~2.1]{Filomat-36-73-1.tex} leads to
\begin{equation}\label{reciprocal=log=der=eq}
\biggl(\frac1{\ln x}\biggr)^{(n)}=\frac1{x^n}\sum_{k=1}^n (-1)^kk!s(n,k)\biggl(\frac1{\ln x}\biggr)^{k+1}, \quad n\in\mathbb{N}.
\end{equation}
This recovers the first formula in~\cite[Lemma~2]{Liu-Qi-Ding-2010-JIS}.
\par
By the way, the formulas~(3.4) and~(3.5) in~\cite[Corollary~3.1]{Filomat-36-73-1.tex} recover the second formula in~\cite[Lemma~2]{Liu-Qi-Ding-2010-JIS}.
\end{rem}

\begin{rem}
In~\cite[Remark~2.2]{Filomat-36-73-1.tex}, it was conjectured that the sequence $a_{n,k}$ for $n\in\mathbb{N}$ and $2\le k\le n+1$ is increasing with respect to $n$ while it is unimodal with respect to $k$ for given $n\ge4$. This conjecture may be partially confirmed as follows.
\par
From~\eqref{a=n=i=eq}, the increasing monotonicity of the sequence $a_{n,k}$ with respect to $n$ follows straightforwardly.
\par
It is clear that the sequence $(k-1)!$ is increasing with $k$ and the sequence
\begin{equation*}
\sum_{\ell_1=1}^{n-1} \frac1{\ell_1}\sum_{\ell_2=1}^{\ell_1-1}\frac1{\ell_2}\dotsm \sum_{\ell_{k-3}=1}^{\ell_{k-4}-1}\frac1{\ell_{k-3}} \sum_{\ell_{k-2}=1}^{\ell_{k-3}-1}\frac1{\ell_{k-2}}
\end{equation*}
is decreasing with $k$. Since $a_{n,n+1}=n!$, see the equation~\eqref{a(n,2)-eq} or~\cite[(2.8)]{Filomat-36-73-1.tex}, we obtain that
\begin{equation}\label{a(n-2)<a(n=n+1)}
a_{n,2}<a_{n,n+1},\quad n\ge2.
\end{equation}
\par
In~\cite[Theorem~2.1]{1st-Sirling-Number-2012.tex}, the integral representation
\begin{equation}\label{n-times-diriv}
s(n,k)=\binom{n}{k}\lim_{x\to0}\frac{\td^{n-k}}{\td x^{n-k}} \biggl\{\biggl[\int_0^\infty\biggl(\int_{1/e}^1 t^{xu-1}\td t\biggr)e^{-u}\td u\biggr]^k\biggr\}
\end{equation}
was created for $1\le k\le n$. Hence,
\begin{align*}
s(n,n-1)&=n\lim_{x\to0}\frac{\td}{\td x} \biggl\{\biggl[\int_0^\infty\biggl(\int_{1/e}^1 t^{xu-1}\td t\biggr)e^{-u}\td u\biggr]^{n-1}\biggr\}\\
&=n(n-1)\lim_{x\to0}\biggl[\int_0^\infty\biggl(\int_{1/e}^1 t^{xu-1}\td t\biggr)e^{-u}\td u\biggr]^{n-2} \\
&\quad\times\lim_{x\to0}\biggl[\int_0^\infty\biggl(\int_{1/e}^1 t^{xu-1}\ln t\td t\biggr)ue^{-u}\td u\biggr]\\
&=n(n-1)\biggl[\int_0^\infty\biggl(\int_{1/e}^1 \frac1t\td t\biggr)e^{-u}\td u\biggr]^{n-2} \int_0^\infty\biggl(\int_{1/e}^1 \frac{\ln t}t\td t\biggr)ue^{-u}\td u\\
&=-\frac{1}{2}n(n-1).
\end{align*}
As a result, by~\eqref{a-stirl-relation}, it follows that
\begin{equation}
a_{n,n}=-(n-1)!s(n,n-1)=\frac{n-1}2n!\ge a_{n,n+1},\quad n\ge3.
\end{equation}
Combining this with~\eqref{a(n-2)<a(n=n+1)} shows that the sequence $a_{n,k}$ for given $n\ge4$ has at leat one maximum with respect to $2<k<n+1$.
\end{rem}

\begin{rem}
By the integral repreaentation~\eqref{n-times-diriv} and direct computation, we can recover that
\begin{align*}
s(n,1)&=\binom{n}{1}\lim_{x\to0}\frac{\td^{n-1}}{\td x^{n-1}} \int_0^\infty\biggl(\int_{1/e}^1 t^{xu-1}\td t\biggr)e^{-u}\td u\\
&=n\lim_{x\to0}\int_0^\infty\biggl[\int_{1/e}^1 t^{xu-1}(\ln t)^{n-1}\td t\biggr]u^{n-1}e^{-u}\td u\\
&=n\int_0^\infty\biggl[\int_{1/e}^1 \frac{(\ln t)^{n-1}}t\td t\biggr]u^{n-1}e^{-u}\td u\\
&=(-1)^{n+1}\int_0^\infty u^{n-1}e^{-u}\td u\\
&=(-1)^{n+1}(n-1)!
\end{align*}
and
\begin{align*}
s(n,2)&=\binom{n}{2}\lim_{x\to0}\frac{\td^{n-2}}{\td x^{n-2}} \biggl\{\biggl[\int_0^\infty\biggl(\int_{1/e}^1 t^{xu-1}\td t\biggr)e^{-u}\td u\biggr]^2\biggr\}\\
&=\binom{n}{2}\lim_{x\to0}\sum_{k=0}^{n-2}\binom{n-2}{k}\int_0^\infty\biggl[\int_{1/e}^1 t^{xu-1}(\ln t)^k\td t\biggr]u^ke^{-u}\td u\\
&\quad\times\int_0^\infty \biggl[\int_{1/e}^1 t^{xu-1}(\ln t)^{n-k-2}\td t\biggr]u^{n-k-2}e^{-u}\td u\\
&=\binom{n}{2}\sum_{k=0}^{n-2}\binom{n-2}{k}\int_0^\infty\biggl[\int_{1/e}^1 \frac{(\ln t)^k}t\td t\biggr]u^ke^{-u}\td u\\
&\quad\times\int_0^\infty \biggl[\int_{1/e}^1 \frac{(\ln t)^{n-k-2}}t\td t\biggr]u^{n-k-2}e^{-u}\td u\\
&=(-1)^n\binom{n}{2}\sum_{k=0}^{n-2}\binom{n-2}{k}\frac{k!}{k+1}\frac{(n-k-2)!}{n-k-1}\\
&=(-1)^n(n-2)!\binom{n}{2}\sum_{k=0}^{n-2}\frac1{(k+1)(n-k-1)}\\
&=(-1)^n\frac{n!}2\sum_{k=0}^{n-2}\frac1{(k+1)(n-k-1)}\\
&=(-1)^n\frac{(n-1)!}2\sum_{k=0}^{n-2}\biggl(\frac1{k+1}+\frac1{n-k-1}\biggr)\\
&=(-1)^n(n-1)!H(n-1),
\end{align*}
where
\begin{equation}
H(n)=\sum_{k=1}^n\frac1n
\end{equation}
is the $n$-th harmonic number. Consequently, we find a relation
\begin{equation}\label{s(n-k)-H(n)-relation}
s(n,2)=(-1)^n(n-1)!H(n-1),\quad n\in\mathbb{N}
\end{equation}
or, equivalently,
\begin{equation}\label{H(n)-s(n-k)-relation}
H(n)=\frac{(-1)^{n+1}s(n+1,2)}{n!},\quad n\in\mathbb{N}
\end{equation}
between the $n$-th harmonic number $H(n)$ and Stirling numbers of the first kind $s(n,2)$.
\par
The relations~\eqref{s(n-k)-H(n)-relation} and~\eqref{H(n)-s(n-k)-relation} may also be deduced by considering~\eqref{a=n=i=eq} and~\eqref{a-stirl-relation} and may also be found in~\cite[p.~275, (6.58)]{GKP-Concrete-Math-2nd}.
\par
For more information on the $n$-th harmonic numbers $H(n)$, please refer to~\cite{harmonic-number-refine-chen.tex} and closely related references therein.
\end{rem}

\begin{rem}
For more information on the second Stirling numbers and the first kind Bernoulli numbers, please refer to~\cite{exp-derivative-sum-Combined.tex, 1st-Sirling-Number-2012.tex} and closely related references therein.
\end{rem}
%


\begin{thebibliography}{9}

\bibitem{GKP-Concrete-Math-2nd}
R. L. Graham, D. E. Knuth, and O. Patashnik, \emph{Concrete Mathematics---A Foundation for Computer Science}, 2nd ed., Addison-Wesley Publishing Company, Reading, MA, 1994.

\bibitem{harmonic-number-refine-chen.tex}
B.-N. Guo and F. Qi, \textit{Sharp bounds for harmonic numbers}, Appl. Math. Comput. \textbf{218} (2011), no.~3, 991\nobreakdash--995; Available online at \url{http://dx.doi.org/10.1016/j.amc.2011.01.089}.

\bibitem{exp-derivative-sum-Combined.tex}
B.-N. Guo and F. Qi, \textit{Some identities and an explicit formula for Bernoulli and Stirling numbers}, J. Comput. Appl. Math. \textbf{255} (2014), 568\nobreakdash--579; Available online at \url{http://dx.doi.org/10.1016/j.cam.2013.06.020}.

\bibitem{Liu-Qi-Ding-2010-JIS}
H.-M. Liu, S.-H. Qi, and S.-Y. Ding, \textit{Some recurrence relations for Cauchy numbers of the first kind}, J. Integer Seq. \textbf{13} (2010), Article~10.3.8.

\bibitem{Nemes-JIS-2011}
G. Nemes, \textit{An asymptotic expansion for the Bernoulli numbers of the second kind}, J. Integer Seq. \textbf{14} (2011), Article~11.4.8.
%

\bibitem{Filomat-36-73-1.tex}
F. Qi, \textit{Explicit formulas for computing Bernoulli numbers of the second kind and Stirling numbers of the first kind}, Filomat \textbf{28} (2014), in press; Available online at \url{http://dx.doi.org/10.2298/FIL??????Q}.

\bibitem{2rd-Bernoulli-2012.tex}
F. Qi, \textit{Explicit formulas for computing Bernoulli numbers of the second kind and Stirling numbers of the first kind}, available online at \url{http://arxiv.org/abs/1301.6845}.

\bibitem{1st-Sirling-Number-2012.tex}
F. Qi, \textit{Integral representations and properties of Stirling numbers of the first kind}, J. Number Theory \textbf{133} (2013), no.~7, 2307\nobreakdash--2319; Available online at \url{http://dx.doi.org/10.1016/j.jnt.2012.12.015}.

\end{thebibliography}
\end{document}